\documentclass[12pt, 14paper,reqno]{amsart}
\usepackage{amsmath,amsfonts,amssymb}
\usepackage[breaklinks]{hyperref}
\usepackage{graphicx}
\usepackage{longtable}
\makeatletter
\@namedef{subjclassname@2020}{%
  \textup{2020} Mathematics Subject Classification}
\makeatother

        \makeatletter
\def\legendre@dash#1#2{\hb@xt@#1{%
  \kern-#2\p@
  \cleaders\hbox{\kern.5\p@
    \vrule\@height.2\p@\@depth.2\p@\@width\p@
    \kern.5\p@}\hfil
  \kern-#2\p@
  }}
\def\@legendre#1#2#3#4#5{\mathopen{}\left(
  \sbox\z@{$\genfrac{}{}{0pt}{#1}{#3#4}{#3#5}$}%
  \dimen@=\wd\z@
  \kern-\p@\vcenter{\box0}\kern-\dimen@\vcenter{\legendre@dash\dimen@{#2}}\kern-\p@
  \right)\mathclose{}}
\newcommand\legendre[2]{\mathchoice
  {\@legendre{0}{1}{}{#1}{#2}}
  {\@legendre{1}{.5}{\vphantom{1}}{#1}{#2}}
  {\@legendre{2}{0}{\vphantom{1}}{#1}{#2}}
  {\@legendre{3}{0}{\vphantom{1}}{#1}{#2}}
}
\def\dlegendre{\@legendre{0}{1}{}}
\def\tlegendre{\@legendre{1}{0.5}{\vphantom{1}}}
\makeatother
\theoremstyle{plain}
\theoremstyle{plain}

\newtheorem{obs}{Observation}

\newtheorem{rem}{Remark}[section]
\numberwithin{equation}{section}
\newtheorem{thm}{Theorem}[section]
\newtheorem{lem}{Lemma}[section]
\newtheorem{cor}{Corollary}[section]
\newtheorem{prop}{Proposition}[section]

\numberwithin{equation}{section}

\usepackage[breaklinks]{hyperref}
\makeatletter
\let\oldHyPsd@CatcodeWarning\HyPsd@CatcodeWarning
\renewcommand{\HyPsd@CatcodeWarning}[1]{
  \ifnum\pdfstrcmp{#1}{math shift}=0    
  \else                                 
    \oldHyPsd@CatcodeWarning{#1}
  \fi
}
\makeatother

\begin{document} 
\title [Elliptic Curve]{Rank of the family of elliptic curves $y^2 = x^3- 5px$}
\author{Arkabrata Ghosh}
\address{G. Arkabrata@ SRM University-AP, Mangalgiri Mandal, Neerukonda, Amravati, Andhra Pradesh-522502.}
\email{arka2686@gmail.com}





\maketitle

\section*{Abstract}

This article considers the family of elliptic curves given by $E_{p}: y^2=x^3-5px$ and certain conditions on an odd prime $p$. More specifically, we have shown that if $p \equiv 7, 23 \pmod {40}$, then the rank of $E_{p}$ is zero for both $ \mathbb{Q} $ and $ \mathbb{Q}(i) $. Furthermore, if the prime $ p $ is of the form $ 40k_1 + 3 $ or $ 40k_2 + 27$ 
where $k_1, k_2 \in \mathbb{Z}$ such that $(5k_1+1)$ or $(5k_2 +4)$ are perfect squares, then the given family of elliptic curves has rank one over $\mathbb{Q}$ and rank two over $\mathbb{Q}(i)$. Moreover, if the prime $ p $ is of the form $ 40k_3 + 11 $ or $ 40k_4 + 19$  where $k_3 ~\text{and}~ k_4 \in \mathbb{Z}$ such that $(160k_3+49)$ or $(160k_4 + 81) $ are perfect squares, then the given family of elliptic curves has rank at least one over $\mathbb{Q}$ and rank at least two over $\mathbb{Q}(i)$. The novelty of this manuscript lies in the fact that we have not used the BSD or parity conjecture to prove the following results in this manuscript.


\section{Introduction}

The arithmetic of the elliptic curve is one of the most fascinating branches of mathematics as it connects number theory to algebraic geometry. 
Recently, there has been a surge among mathematicians in studying the families of elliptic curves given by $E: y^2=  x^3 + bx $ where $ b \in \mathbb{Z}$ and several people are now working in this direction. Chahal \cite{Js88}, in 1988, discussed about Mordell-Weil rank of elliptic curves in detail in his book~(Chapter 7). In 2007, Spearman \cite{SP07} calculated the values of the prime $p$ for which the elliptic curve $ E: y^2 =x^3 -px$ has rank two.  In the same year, he \cite{Sp07} also discovered the condition on $2p$ for which the elliptic curve $y^2= x^3- 2px$ has rank three. In 2010, Hollier et al. \cite{HSY10} considered the elliptic curve of the
form $ E: y^2 = x^3 + pqx$ where $p, q$ are distinct primes. Fujita and Terai \cite{FT11}, in 2011,  considered the elliptic curve of the form $y^2 = x^3- p^{k}x $. where $p$ is a prime and $ k = 1, 2, 3$, and found a necessary and sufficient condition for the rank of the given elliptic curve to be equal to one or two. Daghieh and Didari \cite{DD14} determined that the rank of the elliptic curves of the form $ y^2 = x^3 -3px$, $ p $ is a prime number. They also studied the \cite{DD15} elliptic curves of $y^2 = x^3 -pqx$, where $p$ and $q$ are distinct odd primes. Around the same time, Kim \cite{K15} studied the elliptic curve $y^2 =x^3 \pm 4px$, where $p$ is a prime number. Recently, Mina and Bacani \cite{MB23} studied the curve $y^2= x^3 -3pqx$ for distinct odd primes $p$ and $q$ and found two different sufficient conditions on $p$ and $q$ such that the given elliptic curve, under those conditions, has rank zero and one, respectively. Finally, Ghosh \cite{Gh25} studied the torsion and rank of the family of elliptic curves given by $ y^2 = x^3 - 5pqx$.

In this article, we aim to find the rank of the following family of elliptic curves:
\begin{equation}
    \label{eq:1.1}
    E_{p}: y^2 = x^3- 5px,
\end{equation}
where $p$ is an odd prime that satisfies certain conditions. First, we will provide conditions for $p$ such that the rank of the family elliptic curve $E_{p}$ given by equation \eqref{eq:1.1} is zero. More precisely, we will prove the following Theorem.

\begin{thm}
    \label{Thm:1.1}
    
    If $p$ is an odd prime and $ p $ is congruent either to  $7~\text{or to}~ 23 \pmod {40} $, then the rank of the elliptic curve given by \eqref{eq:1.1} is zero over $\mathbb{Q}$.
    
\end{thm}

Now we will update the conditions of $p$  so that the rank of the elliptic curve given by \eqref{eq:1.1} is exactly one. This is stated in the following Theorem.

\begin{thm}
    \label{Thm:1.2}

 Suppose $p$ is an odd prime. If $ p = 40k_1 +3$ where $ (5k_1+1)$ is a perfect square or $ p = 40k_2 +27$ where $ (5k_2+4) $ is a perfect square, then the rank of the elliptic curve given by \eqref{eq:1.1} is one over $\mathbb{Q}$.
 \end{thm}

 We finally end this manuscript with the following result.

 \begin{thm}
     \label{Thm:1.3}
 
 Suppose $p$ is an odd prime. If $ p=40k_3+11$ where $(160k_3 +49)$ is a perfect square or  $p=40k_4+19$ where $(160k_4 +81)$ is a perfect square, then the rank of the elliptic curve given by \eqref{eq:1.1} is at least one over $\mathbb{Q}$.

 \end{thm}

 In \S 2, we shall discuss some preliminaries regarding the Mordeil-Weil rank of an elliptic curve, and finally in \S 3, we shall prove Theorems \ref{Thm:1.1}, \ref{Thm:1.2}, and \ref{Thm:1.3}


   

\section{preliminaries}

Let $E$ be an elliptic curve over the field $\mathbb{K}$ of characteristic not equal to $2$ or $3$ and let $ E(\mathbb{K})$ denote the $\mathbb{K}$-rational points of $E$ over $\mathbb{K}$. Mordell-Weil Theorem asserts that $  E(\mathbb{K})$ is a finitely generated abelian group and can be represented as

\begin{equation*}
    E(\mathbb{K}) \cong \mathbb{Z}^r \oplus E(\mathbb{K})_{tors}, 
\end{equation*}
where $ r \geq 0$ is called the rank of the elliptic curve $E$ over $\mathbb{K}$ and $ E(\mathbb{K})_{tors}$ denotes the torsion subgroup of $ E(\mathbb{K})$, which is a finite abelian group consisting of elements of finite order. Mazur's Theorem \cite{M77} gives an explicit idea about the torsion subgroup $ E(\mathbb{K})_{tors}$ where $\mathbb{K}=\mathbb{Q}$.





Now to compute the rank of the curve given by \eqref{eq:1.1}, we need to use the method of $2$-descent. We will describe it briefly here, and one can look at (see \cite{ST15}) for more details. Suppose that $E : Y^2= X^3 + a X^2 + bX$ is an elliptic curve over $\mathbb{Q}$ and $\overline{E}: Y^2 = X^3 - 2aX^2 + (a^2- 4b) X$ is the corresponding isogenous curve to $E$. Hence, there exists an isogeny $ \phi: E \rightarrow \overline{E}$ of degree $2$ given by
\begin{equation*}
    \phi(x,y) ~ = ~ \bigg(\frac{y^2}{x^2}, \frac{y(x^2-b)}{x^2}\bigg).
\end{equation*}
Moreover, let $\mathbb{Q}^{*}$ be the multiplicative group of all non-zero rational numbers, and $\mathbb{Q^*}^{2}$ be its subgroup of squares of elements of $ \mathbb{Q^*}$. Hence, $ \mathbb{Q^*}/\mathbb{Q^*}^{2}$ is the multiplicative group of all non-zero rational numbers modulo squares. Also, we denote set of rational points on $E$ and $\overline{E}$ by $\Gamma$ an $ \overline{\Gamma}$ respectively. Now we define the $2$-descent homomorphism $\alpha: \Gamma \rightarrow \mathbb{Q^*}/\mathbb{Q^*}^{2}$ by

\begin{equation*}
    \alpha(P)= 
               \begin{cases}
                 1 \pmod {\mathbb{Q^*}^2}, ~\text{if}~ P= \mathcal{O}, ~\text{the point at infinity},\\
                 b \pmod {\mathbb{Q^*}^2}, ~\text{if}~ P= (0,0),\\
                  x \pmod {\mathbb{Q^*}^2}, ~\text{if}~ P= (x,y ) ~\text{with}~ x \neq 0.
                 \end{cases}
\end{equation*}

Similarly, we can define the $2$-decent homomorphism on the isogeneous curve $\overline{E}(\mathbb{Q}) $ as follows: $\overline{\alpha}:  \overline{\Gamma} \rightarrow  \mathbb{Q^*}/\mathbb{Q^*}^{2}$ by 

\begin{equation*}
    \overline{\alpha}(\overline{P})= 
               \begin{cases}
                 1 \pmod {\mathbb{Q^*}^2}, ~\text{if}~ \overline{P}= \mathcal{\overline{O}}, ~\text{the point at infinity},\\
                 \overline{b} \pmod {\mathbb{Q^*}^2}, ~\text{if}~ \overline{P}= (0,0),\\
                 x \pmod {\mathbb{Q^*}^2}, ~\text{if}~ \overline{P}= (x,y ) ~\text{with}~ x \neq 0,
                 \end{cases}
\end{equation*}

where $\overline{b}= a^2 - 4b$.

The group $\alpha(\Gamma)$ consists of $1,b$ and all factors $b_1$ of $b$, all modulo $\mathbb{Q^*}^2$. Here $b_1 \neq 1, ~\text{or}~ b \pmod {\mathbb{Q^*}^2}$, such that the triple $(N, M, e) \in \mathbb{Z}^3$, where $ M \neq 0, e \neq 0$ solves the Diophantine equation ( or 'torsors')

\begin{equation*}
    \mathcal{T}: N^2 = b_1 M^4 + a M^2 e^2 + b_2 e^4, ~\text{with}~  b_1 b_2 =b,
\end{equation*}

and satisfies the following criterion :

\begin{equation*}
    \gcd(N,e)=\gcd(M,e)=\gcd(b_1, e) =\gcd(b_2, M)= \gcd(M,N)=1.
\end{equation*}

Similarly, the group $\overline{\alpha}(\overline{\Gamma})$ consists of $1, a^2-4b$ and all factors $b_1$ of $a^2-4b$, all modulo $\mathbb{Q^*}^2$. Here $b_1 \neq 1, ~\text{or}~ a^2-4b \pmod {\mathbb{Q^*}^2}$, such that the triple $(N, M, e) \in \mathbb{Z}^3$, where $ M \neq 0, e \neq 0$ solves the Diophantine equation ( or 'torsors')

\begin{equation*}
    \mathcal{T'}: N^2 = b_1 M^4 -2a M^2 e^2 + b_2 e^4, ~\text{with}~  b_1 b_2 =a^2-4b,
\end{equation*}

and the same GCD criterion mentioned above. Now to compute the rank of an elliptic curve $E$, we use the following Proposition(see \cite{ST15}).
\begin{prop}
    \label{prop:2.1} Let $r$ be the rank of $E(\mathbb{Q}) $ and $\alpha$ and $\Bar{\alpha}$ are as above. Then
    $$
       \frac{1}{4} |\alpha(\Gamma)| \cdot |\overline{\alpha}(\overline{\Gamma})|= 2^r.
    $$
\end{prop}

Now, if $\mathbb{K} = \mathbb{Q}(\sqrt{m})$, where $m$ is a square-free integer, we can find the rank of any elliptic curve $E$ over $\mathbb{K}$ by adding ranks of $E$ and its' $m$-twist $E[m]$ over $\mathbb{Q}$. This can be seen in the following result (see \cite{S09}).

\begin{thm}
    \label{thm:2.1}
    Let $\mathbb{K} = \mathbb{Q}(\sqrt{m})$ be a quadratic field, where $m$ is a square-free integer. Let $E: y^2 = x^3 + ax^2 + bx $ be an elliptic curve over $\mathbb{Q}$ and $E[m]: y^2 = x^3 + max^2 + m^2 bx$ be the $m$-twist of $E$. Then

    $$
    rank(E(\mathbb{K}))= rank(E(\mathbb{Q})) + rank(E[m](\mathbb{Q}) ).
    $$
    
\end{thm}

We shall end this section with the following corollary.

\begin{cor}
 \label{cor:3.2}
 Consider the field $\mathbb{K}=\mathbb{Q}(i)$. Then we can say the following about the rank of the elliptic curve $E_{p}$ given by \ref{eq:1.1} over $\mathbb{K}$. Then the following holds.

 (i) If $ p$ is an odd prime that meets the condition as in Theorem \ref{Thm:1.1}, then the rank of $E_{p} $ over $\mathbb{K}$ is zero.

 (ii) If $ p$ is an odd prime and and $k_1, k_2 \in \mathbb{Z}$  such that they satisfy the hypothesis of Theorem \ref{Thm:1.2}, then the rank of $E_{p} $ over $\mathbb{K}$ is exactly $2$.

 (iii) If $ p$ is an odd prime and and $k_3 \in \mathbb{Z}$  such that they satisfy the hypothesis of Theorem \ref{Thm:1.3}, then the rank of $E_{p} $ over $\mathbb{K}$ is at least $2$.

\end{cor}
    
For some applications of elliptic curves, we refer the reader to \cite{IM22} and \cite{SS15}.

\section{Rank of $\texorpdfstring{E_{p}}{}$}

We start this section with the following remark.

\begin{rem}
\label{rem:3.1}
Let $p$ be an odd prime. If $ p \equiv 7, 23 \pmod {40}$, then $ \genfrac(){}{0}{p}{5}=-1$.
\end{rem}

\begin{proof}
    From the definition of the Legendre symbol, we know
    \begin{equation*}
       \genfrac(){}{0}{p}{5} \genfrac(){}{0}{5}{p} =(-1)^{\frac{p-1}{2} \frac{5-1}{2}}= (-1) ^{(p-1)} =1.
    \end{equation*}
   As $p \equiv 7,23 \pmod {40} \equiv 2,3 \pmod 5$,  $ \genfrac(){}{0}{5}{p}=-1$. Then from the above equation we can say $ \genfrac(){}{0}{p}{5} = -1$.
\end{proof}


From equation \eqref{eq:1.1}, we have $  E_{p}: ~ y^2= x^3-5px $ and the corresponding $\overline{E}_{p}: y^2 = x^3 + 20px$. To compute the rank of $E_{p}$, we first need to determine $|\alpha(\Gamma)|$. Note that $ 1, -5p  \in \alpha(\Gamma)$ by the definition of $\alpha$. In that case, the following set gives all possible divisors $b_1$ of $-5p$ other than $1 ~\text{and}~ -5p$  modulo $\mathbb{Q^*}^2$.

\begin{equation*}
    S= \{-1, \pm 5, \pm p, 5p \}.
\end{equation*}

We then consider the solvability of the following torsors over a set of integers.

\begin{equation*}
  \begin{split}
    \mathcal{T}_1: & N^2 = 5pM^4 - e^4 \\
    \mathcal{T}_2 : & N^2 = 5M^4- pe^4 \\
    \mathcal{T}_3 : & N^2 = pM^4 - 5e^4\\
    \end{split}
\end{equation*}

\begin{lem}
    \label{le:3.1}
    For $ p \equiv 7, 23 \pmod{40}$, there are no integer solutions for the torsor $ \mathcal{T}_1: N^2=5pM^4-e^4$,
\end{lem}

\begin{proof}
     Reducing $ \mathcal{T}_1 $ modulo $p$, we get $N^2 \equiv -e^4 \pmod{p}$. It implies 

    \begin{equation*}
        1 = \genfrac(){}{0}{-e^4}{p}= \genfrac(){}{0}{-1}{p} \genfrac(){}{0}{e^4}{p}= \genfrac(){}{0}{-1}{p}
    \end{equation*}

Now $ \genfrac(){}{0}{-1}{p}=1 \iff p \equiv 1 \pmod 4$. As $ p \equiv 7, 23 \pmod {40} \equiv 3 \pmod 4$, we arrived at a contradiction. Hence, we can say that $  \mathcal{T}_1  $ has no solution in $ \mathbb{Z}$.
    
\end{proof}

\begin{lem}
\label{le:3.2}
For $p \equiv 7, 23 \pmod{40}$, there are no integer solutions for the torsors $ \mathcal{T}_2: N^2=5M^4-pe^4$.
\end{lem}

\begin{proof}
    Reducing $ \mathcal{T}_2$ modulo $5$, we get $N^2 \equiv -pe^4 \pmod 5$. So,

    \begin{equation*}
        1 = \genfrac(){}{0}{-pe^4}{5}= \genfrac(){}{0}{-1}{5} \genfrac(){}{0}{p}{5} \genfrac(){}{0}{e^4}{5} =  \genfrac(){}{0}{p}{5}.
    \end{equation*}
    This is a contradiction, as from the remark \eqref{rem:3.1}, 
we know $\genfrac(){}{0}{p}{5}=-1$.
\end{proof}

\begin{lem}
\label{le:3.3}
For $p \equiv 7, 23 \pmod{40}$, there are no integer solutions for the torsors $ \mathcal{T}_3: N^2=pM^4-5e^4$.
\end{lem}

\begin{proof}
    Reducing $ \mathcal{T}_3$ modulo $5$, we get $N^2 \equiv pM^4 \pmod 5$. So,

    \begin{equation*}
        1 = \genfrac(){}{0}{pM^4}{5}=  \genfrac(){}{0}{p}{5} \genfrac(){}{0}{M^4}{5} =  \genfrac(){}{0}{p}{5}.
    \end{equation*}
    This is a contradiction as from the remark \eqref{rem:3.1}, 
we know $\genfrac(){}{0}{p}{5}=-1$.

\end{proof}

Hence, combining lemmas \ref{le:3.1}-\ref{le:3.3}, we get $ \alpha(\Gamma)= \{1, -5p \}$ and hence $ |\alpha(\Gamma)| = 2$. 

Now we move on to find $ |\overline{\alpha}(\overline{\Gamma})|$. As $\overline{E}_{p}: y^2= x^3 + 20px$, we know that $ \{1, 5p \} \in \overline{\alpha}(\overline{\Gamma})$. Now, all possible divisors $b_1$ of $20p$ except $1$ and $5p$ modulo $ \mathbb{Q^*}^2$ are elements of the following set.

\begin{equation*}
    T=\{2, 5 ,10, p, 2p, 10p\}.
\end{equation*}

We deliberately remove any negative values of $b_
1$ from the set $T$ as the corresponding torsors $N^2= b_1 M^4 + b_2 e^4$ will have no solutions if both $b_1$ and $b_2$ are negative. Hence, we will take all possible values of $b_1$ from the set $T$ and consider the solvability of the corresponding torsor over $\mathbb{Z}$. They are as follows:

\begin{equation*}
    \begin{split}
        \mathcal{T}_1': & N^2 = 2M^4+10pe^4 \\
        \mathcal{T}_2': & N^2 = 20M^4+ pe^4 \\
        \mathcal{T}_3': & N^2 = 2pM^4+10e^4 \\
        \mathcal{T}_4': & N^2= 5M^4 +  4pe^4 
    \end{split}
\end{equation*}

\begin{lem}
    \label{le:3.4}
    For $p \equiv 7, 23 \pmod{40}$, there are no integer solutions for the torsor $ \mathcal{T}_1': N^2=2M^4+10pe^4$.
\end{lem}
\begin{proof}
    Reducing $ \mathcal{T}_1' $ modulo $5$, we get $N^2 \equiv 2M^4 \pmod 5$. Hence,

    \begin{equation*}
        1= \genfrac(){}{0}{2M^4}{5}= \genfrac(){}{0}{2}{5} \genfrac(){}{0}{M^4}{5}=  \genfrac(){}{0}{2}{5}=-1, 
    \end{equation*}
a contradiction. Thus,  $ \mathcal{T}_1' $ has no solution in $\mathbb{Z}$.
\end{proof}

\begin{lem}
    \label{le:3.6} 
    For $p \equiv 7, 23 \pmod{40}$, there are no integer solutions for the torsors $ \mathcal{T}_2': N^2=20M^4+ pe^4$ and $ \mathcal{T}_4': N^2=5M^4+ 4pe^4$.

\end{lem}

\begin{proof}
    Taking into account $ \mathcal{T}_2'$ modulo $p$, we get $N^2 \equiv 5M^4 \pmod p$. Hence,

    \begin{equation*}
        1= \genfrac(){}{0}{20M^4}{p}= \genfrac(){}{0}{5}{p}.
    \end{equation*}
    As $ p \equiv 7, 23 \pmod {40} \equiv 2, 3 \pmod 5$, we have $ \genfrac(){}{0}{5}{p}=-1 $. So we arrive at a contradiction and can conclude that $\mathcal{T}_2'$ has no solution in $\mathbb{Z}$.

    Similarly, reducing $\mathcal{T}_4'$ modulo $p$, we get $\genfrac(){}{0}{5}{p}=1 $, a contradiction. Hence, we can say that $\mathcal{T}_4'$ has no solution in integers.
\end{proof}    

\begin{lem}
    \label{le:3.7} 
     For $p \equiv 7, 23 \pmod{40}$, there are no integer solutions for the torsor $ \mathcal{T}_3': N^2=2pM^4+10e^4$.

\end{lem}
\begin{proof}
    As $N^2=2(pM^4 + 5e^4)$, we know $N$ is even and as $\gcd(N,M)=\gcd(N,e)=1$. Hence, both $M$ and $e$ are odd. Assuming $N=2N_1$ for some $N_1 \in \mathbb{Z}$, we get $2N_1^2= pM^4 + 5e^4$. As both $M$ and $e$ are odd, we get $M^4 \equiv e^4 \equiv 1 \pmod 8$ and using it, we can say $ 2N_1^2 \equiv p + 5 \pmod 8$. By our assumption, $ p \equiv 7, 23 \pmod {40} \equiv 7 \pmod 8$. Hence, $2N_1^2 \equiv  4 \pmod 8$ which implies $N_1^2 \equiv 2 \pmod 4$, a contradiction. Thus, $\mathcal{T}_3'$ has no solutions in integers.
\end{proof}

So, from the Lemmas \ref{le:3.4}-\ref{le:3.7}, we can conclude that $ \overline{\alpha}(\overline{\Gamma}) = \{1, 5p \}$ and hence $ |\overline{\alpha}(\overline{\Gamma})|=2$. 

We are ready to prove Theorem \ref{Thm:1.1}, which is as follows.

\begin{proof}
    As $ |\alpha(\Gamma)|= 2=|\overline{\alpha}(\overline{\Gamma})| $,  from the Proposition \ref{prop:2.1}, we have

    \begin{equation}
        2^{r}= \frac{1}{4}( |\alpha(\Gamma)| \cdot |\overline{\alpha}(\overline{\Gamma})|=1,
    \end{equation}
    from which we get $r=r(E_{p})=0$.
\end{proof}

The following table \ref{tab:1} confirms the result of Theorem \eqref{Thm:1.1}. Here we take some primes $p$  which satisfy the condition of Theorem \ref{Thm:1.1} and list down the rank of the corresponding elliptic curve $E_{p}$. All computations are done using SAGE \cite{SA}.


\begin{table}[ht]
    \centering

\begin{tabular}{|c| c|} 

\hline
 $p$ & $\text{rank of}~ E_{p}$ \\
\hline 
 7 & 0 \\
\hline
 47 & 0 \\
\hline
 23 & 0 \\
\hline
 103 & 0 \\
\hline 

\end{tabular}
 \caption{The values of $p$ that verifies the condition of Theorem \ref{Thm:1.1}}
    \label{tab:1}
\end{table}


\begin{cor}
    \label{cor:3.1} Under the same assumption as of Theorem \ref{Thm:1.1}, if the congruence $p \equiv 7,23 \pmod {40}$, is replaced by $ p \equiv 27,3  \pmod {40}$, then the rank of elliptic curve $E_{p}$ given by equation \eqref{eq:1.1} is at most one. 
\end{cor}

\begin{proof}
    A simple calculation shows that $ p \equiv 27 \pmod {40} \equiv 3 \pmod {8}$ and $p \equiv 3 \pmod {40} \equiv 3 \pmod 8$. Throughout all lemmas \ref{le:3.1}-\ref{le:3.7}, we have used the condition that $ p \equiv 2 \pmod{5}$, but not  $ p \equiv 7 \pmod 8$ except lemma \ref{le:3.7}. That particular lemma uses the condition $ p \equiv 7 \pmod 8$ and, hence, the torsor $\mathbf{T}_3'$ used the same assumption. According to this new assumption, $p$ does not satisfy $ p \equiv 7 \pmod {8}$ but rather assumes $ p \equiv 3 \pmod{8}$. This implies the torsor $\mathbf{T}_3'$ may have integer solutions under this new condition. So, from lemma \ref{le:3.7}, we can say that it is possible $ 2p ~\text{and}~ 10 \in \overline{\alpha}(\overline{\Gamma})$, making $ |\overline{\alpha}(\overline{\Gamma})| \leq 4$. As a consequence, we obtain $ o \leq r \leq 1$.

\end{proof}

Using this corollary, we will now prove Theorem \ref{Thm:1.2}

\begin{proof}
 If $ p = 40 k_1 + 27$ for some $k_1 \in \mathbb{Z}$, then it satisfies the assumptions of the corollary \ref{cor:3.1}. So, the torsor that may have integer solutions is $ \mathcal{T}_{5}'$. Now 

\begin{equation*}
    \mathcal{T}_{3}': N^2 = 2( 40k_1 +27) + 10= 16(5k_1+4),
\end{equation*}

implies $ \mathcal{T}_{3}'$ has a solution $(N,M,e)= ( 4 \sqrt{5k_1+4}, 1,1)$. Similarly, if $ p = 40 k_2 + 3$ for some $k_2 \in \mathbb{Z}$, then it also satisfies the assumptions of corollary \ref{cor:3.1}. So, the torsor that may have integer solutions  $ \mathcal{T}_{3}'$. Now 

\begin{equation*}
    \mathcal{T}_{3}': N^2 = 2( 40k_2 +3) + 10= 16(5k_2+1),
\end{equation*}

implies $ \mathcal{T}_{3}'$ has a solution $(N,M,e)= ( 4 \sqrt{5k_2+1}, 1,1)$.

Hence, in both cases, $ 2p ~\text{and}~ 10 \in |\overline{\alpha}(\overline{\Gamma})| $. Hence,  $ |\overline{\alpha}(\overline{\Gamma})|=4$ and as a consequence, using the Proposition \ref{prop:2.1}, we obtain $ r(E_{p})=1$.

\end{proof}

The following table \ref{tab:2} confirms the result of Theorem \ref{Thm:1.2}. Here we take some primes $p$  which satisfy the condition of Theorem \ref{Thm:1.2} and list the rank of the corresponding elliptic curve $E_{p}$. All calculations are performed using SAGE \cite{SA}.

\begin{table}[ht]
    \centering
\begin{tabular}{|c| c| c|} 

\hline
$k_1/ k_2$ & $p$ & $\text{rank of}~ E_{p}$ \\
\hline 
 
 $k_1=0$ & 3 & 1 \\
\hline
 $k_2=1$ & 67 & 1 \\
\hline
 $k_1=7$ & 283 & 1 \\
\hline
$k_1=16$ & 643 & 1\\
\hline
$k_2=145$ & 5827 & 1\\
\hline 

\end{tabular}
 \caption{The values of $p$ that verifies conditions of Theorem \ref{Thm:1.2}}
    \label{tab:2}
\end{table}

We shall prove our final Theorem \ref{Thm:1.3}.

\begin{proof}

From the discussion in \S 2, we know $1, -5p \in \alpha(\Gamma)$ and $ 1, 5p  \in \ \overline{\alpha}(\overline{\Gamma})$. If $ p = 40 k_3 + 11$ for some $k_3 \in \mathbb{Z}$, then from Torsor $\mathcal{T}_4'$, we have the following. 
  \begin{equation*}
    \mathcal{T}_{2}': N^2 = 2^4( 40k_3 +11) + 20= 640k_3 +196= 4(160k_3+49).
\end{equation*}
It implies $\mathcal{T}_2'$ has a solution $(N,M,e) =(2 \sqrt{160k_3+49},1,2) $. Hence, we can say that $ 20 ~\text{and}~p \in  \overline{\alpha}(\overline{\Gamma})$. 

Similarly, if $  p = 40 k_4 + 19$ for some $k_4 \in \mathbb{Z}$, then from Torsor $\mathcal{T}_4'$,we have the following. 
  \begin{equation*}
    \mathcal{T}_{2}': N^2 = 2^4( 40k_4 +19) + 20= 640k_4 +324= 4(160k_3+81).
\end{equation*}
It implies $\mathcal{T}_2'$ has a solution $(N,M,e) =(2 \sqrt{160k_4+81},1,2) $. Hence, we can say that $ 20 ~\text{and}~p \in \overline{\alpha}(\overline{\Gamma})$.





%

Hence, in both cases, $ p ~\text{and}~ 5 \in \overline{\alpha}(\overline{\Gamma}) $. Hence, we can say $\{1, -5p \}  \subset  \alpha(\Gamma)$ and $ \{1, 5, p, 5p \} \subset  \overline{\alpha}(\overline{\Gamma})$. Hence, $|\alpha(\Gamma)| \geq 2 $ and $|\overline{\alpha}(\overline{\Gamma})| \geq 4$. Now, using Proposition \ref{prop:2.1}, we obtain $ r(E_{p}) \geq 1$.

\end{proof}

The following table \ref{tab:3} confirms the result of Theorem \ref{Thm:1.3}. Here we take some primes $p$  which satisfy the condition of Theorem \ref{Thm:1.3} and list the rank of the corresponding elliptic curve $E_{p}$. All calculations are performed using SAGE \cite{SA}.


\begin{table}[ht]
\label{tab:3}
    \centering
\begin{tabular}{|c| c| c| c|} 

\hline
$k_3/k_4$ & $N$ & $p$ & $\text{rank of}~E_{p}$ \\
\hline
$k_3= 0$ & 14 & 11 & 1 \\
\hline
 $k_3=3$ & 46 & 131 & 1 \\
\hline 
$k_4=0$ & 18 & 19 & 1 \\
\hline
$k_4=9$ & 78 & 379 & 1\\
\hline
\end{tabular}
 \caption{The values of $p$ that verifies conditions of Theorem \ref{Thm:1.3}}
\end{table}

Now using Theorems \ref{Thm:1.1}, \ref{Thm:1.2}, \ref{Thm:1.3} 
and,  \ref{thm:2.1}, we will prove the corollary \ref{cor:3.2}

\begin{proof}
    (i) If $\mathbb{K}=\mathbb{Q}(i)$, then from Theorem \ref{thm:2.1}, we have $m=-1$ and, hence, $E_{p}[-1]: y^2 =x^3 -(-1)^2 ~5px= x^3 -5px= E_{p}$. Hence, using Theorem \ref{thm:2.1}, we can say 

    \begin{equation*}
        \begin{split}
            rank(E_{p}(\mathbb{K})) &= rank(E_{p}(\mathbb{Q})) + rank(E_{p}[-1](\mathbb{Q}))\\
            & = 2~ rank(E_{p}(\mathbb{Q})).
        \end{split}
    \end{equation*}
From Theorem \ref{Thm:1.1}, we know that $rank(E_{p}(\mathbb{Q}))=0$. So, our claim follows directly from the above equation.

(ii) From Theorem \ref{Thm:1.2}, we know that $ rank(E_{p}(\mathbb{Q}))=1$. So, using the above equation $   rank(\mathbb{K}) = 2~ rank(E_{p}(\mathbb{Q}))$, we get the desired conclusion. 

(iii) From Theorem \ref{Thm:1.3}, we know that $ rank(E_{p}(\mathbb{Q})) \geq 1$. Hence, using the above equation $   rank(\mathbb{K}) = 2~ rank(E_{p}(\mathbb{Q}))$, we get the desired conclusion. 

\end{proof}

\section*{Concluding Remarks}

In this manuscript, we have shown that for different conditions on an odd prime $p$, the rank of the elliptic curve given by \eqref{eq:1.1} is zero, one, or at least one over $\mathbb{Q}$. Now the question remains, can we put some conditions on $p$ such that the rank of the elliptic curve given by \eqref{eq:1.1} is $2$ or higher? In the case of rank two, we can state the following observation without proof, and we encourage readers to prove this.

\begin{obs}
\label{obs:1}
    If $p$ is an odd prime and $ p \equiv 31 \pmod {80}$, then the rank of the elliptic curve given by \eqref{eq:1.1} is two.
\end{obs}

The following table \ref{tab:4} supports the above observation. All computations have been done using SAGE \cite{SA}.

\begin{table}[ht]
\label{tab:4}
    \centering
\begin{tabular}{|c| c|} 

\hline
 $p$ & $\text{rank of}~E_{p}$ \\
\hline
31 & 2 \\
\hline
 191 & 2 \\
\hline 
271 & 2 \\
\hline
431 & 2\\
\hline
\end{tabular}
 \caption{The values of $p$ that satisfies conditions of observation \ref{obs:1}}
\end{table}

\section*{Funding and Conflict of Interests/Competing Interests} The author has no financial or non-financial interests to disclose that are directly or indirectly related to the work. The author has no funding sources to report.

\section*{Data availability statement} No outside data was used to prepare this manuscript.

\section*{Acknowledgment} The author acknowledges Professor Kalyan Chakraborty and Professor Paul Voutier for their continuous support during the preparation of this manuscript. He also acknowledges the support of Dr. Kalyan Banerjee for his constant support during the difficult time, when the author was fighting cancer.

\end{document}